%% file: fokker_planck__arXiv.tex
\documentclass[]{scrartcl}

\usepackage[utf8]{luainputenc}
\usepackage[USenglish]{babel}
\usepackage{csquotes}

\usepackage[a4paper,top=27mm,bottom=20mm,inner=25mm,outer=20mm]{geometry}

\usepackage[%
  backend=bibtex,bibencoding=ascii,
  style=numeric-comp,
  giveninits=true, uniquename=init, 
  natbib=true,
  url=true,
  doi=true,
  isbn=false,
  backref=false,
  maxnames=99,
  ]{biblatex}
\usepackage{amsmath}
\allowdisplaybreaks
\numberwithin{equation}{section}
\usepackage{amssymb}
\usepackage{commath}
\usepackage{mathtools}
\usepackage{bbm}
\usepackage{nicefrac}
\usepackage{subdepth}

\usepackage{siunitx}
\usepackage{algorithm}
\usepackage{algpseudocode}

\usepackage{amsthm}
\usepackage{thmtools}
\usepackage{etoolbox}
\makeatletter
\patchcmd{\thmt@setheadstyle}
 {\bgroup\thmt@space}
 {\thmt@space}
 {}{}
\patchcmd{\thmt@setheadstyle}
 {\egroup\fi}
 {\fi}
 {}{}
\makeatother
\declaretheoremstyle[
  bodyfont=\normalfont\itshape,
  headformat=\NAME\ \NUMBER\NOTE,
]{myplain}
\declaretheoremstyle[
  headformat=\NAME\ \NUMBER\NOTE,
]{mydefinition}
\newcommand{\envqed}{{\lower-0.3ex\hbox{$\triangleleft$}}}
\declaretheorem[style=myplain,numberwithin=section]{theorem}

\declaretheorem[style=myplain,numberlike=theorem]{proposition}

\declaretheorem[style=mydefinition,numberlike=theorem]{definition}

\usepackage[plainpages=false,pdfpagelabels,hidelinks,unicode]{hyperref}

\usepackage{color}
\usepackage{graphicx}
\usepackage[small]{caption}
\usepackage{subcaption}

\begingroup\expandafter\expandafter\expandafter\endgroup
\expandafter\ifx\csname pdfsuppresswarningpagegroup\endcsname\relax
\else
\fi

\usepackage{booktabs}
\usepackage{rotating}
\usepackage{multirow}

\usepackage{enumitem}

\usepackage{ifluatex}
\ifluatex
  \usepackage[no-math]{fontspec}
\else
  \usepackage[T1]{fontenc}
\fi
\usepackage{newpxtext,newpxmath}

\usepackage{bm}
\usepackage{cleveref}
\input{macros}

\makeatletter
\renewcommand*{\@fnsymbol}[1]{\ensuremath{\ifcase#1\or \mathsection\or \mathparagraph\or \|\or *\or **\or \dagger\dagger
   \or \ddagger\ddagger \else\@ctrerr\fi}}
\makeatother

\newcommand{\orcid}[1]{ORCID:~\href{https://orcid.org/#1}{#1}}
\usepackage{authblk}

\newenvironment{keywords}{\par\textbf{Key words.}}{\par}
\newenvironment{AMS}{\par\textbf{AMS subject classification.}}{\par}

\title{Structure-Preserving Numerical Methods for Fokker-Planck Equations}

\author[1]{Hanna Bartel\thanks{\orcid{0009-0001-1889-0123}}}
\affil[1]{Applied Mathematics, University of Hamburg, Germany}

\author[1]{Joshua Lampert\thanks{\orcid{0009-0007-0971-6709}}}

\author[2]{Hendrik~Ranocha\thanks{\orcid{0000-0002-3456-2277}}}
\affil[2]{Institute of Mathematics, Johannes Gutenberg University Mainz, Germany}

\date{November 04, 2024}

\makeatletter
\hypersetup{pdfauthor={Hanna Bartel, Joshua Lampert, Hendrik Ranocha}}
\hypersetup{pdftitle={Structure-Preserving Numerical Methods for Fokker-Planck Equations}}
\makeatother

\begin{document}

\maketitle

\begin{abstract}
\noindent
  \input{abstract.tex}
\end{abstract}



\input{fokker_planck}

\section*{Acknowledgments}

\input{acknowledgments}

\printbibliography

\end{document}

%% file: macros.tex
\usepackage{tikz}
\usepackage{graphicx}
\usepackage{amssymb} 
\usepackage{booktabs}

\newcommand{\fp}{Fokker-Planck }
\newcommand{\fpe}{Fokker-Planck equation }
\newcommand{\fpeo}{Fokker-Planck equation}
\newcommand{\fpes}{Fokker-Planck equations }
\newcommand{\fpeso}{Fokker-Planck equations}

\newcommand{\ccm}{Chang-Cooper method }
\newcommand{\ccmo}{Chang-Cooper method}

\newcommand{\dw}{\Delta_w}

\newcommand{\iph}{_{i+\frac{1}{2}}}
\newcommand{\imh}{_{i-\frac{1}{2}}}

\newcommand{\liph}{\lambda_{i+\frac{1}{2}}}

\newcommand{\pii}{p_{i,j}}

\newcommand{\dii}{d_{i,j}}
\newcommand{\djj}{d_{j,i}}
\newcommand{\sM}{\sum_{i\in\mathcal{M}}}
\newcommand{\sMj}{\sum_{j\in\mathcal{M}}}

\newcommand{\dt}{\Delta_t}
\newcommand{\Mo}{\overline{\mathcal{M}}}

\newcommand{\mpe}{modified Patankar-Euler scheme }

\newcommand{\mpeo}{modified Patankar-Euler scheme}
\newcommand{\mper}{modified Patankar-Euler scheme \eqref{eq:mpe} }
\newcommand{\mprk}{modified Patankar-Runge-Kutta scheme }

\newcommand{\mprko}{modified Patankar-Runge-Kutta scheme}
\newcommand{\mprkr}{modified Patankar-Runge-Kutta scheme \eqref{eq:mprk} }
\newcommand{\ees}{explicit Euler scheme }
\newcommand{\eeso}{explicit Euler scheme}

\newcommand{\tww}{$\dt=\Delta_w^2 / (2\sigma^2)$ }
\newcommand{\twwo}{$\dt=\Delta_w^2 / (2\sigma^2)$}

\newcommand{\mN}{\mathbb{N}}

\newcommand{\mR}{\mathbb{R}}

\newcommand{\mf}{\mathcal{F}}

\newcommand{\dd}{\:\mathrm{d}}

%% file: abstract.tex
A common way to numerically solve Fokker-Planck equations is the Chang-Cooper method 
in space combined with one of the Euler methods in time. However, the explicit Euler method is only conditionally positive, leading to severe restrictions on the time step to ensure positivity.
On the other hand, the implicit Euler method is robust but nonlinearly implicit. Instead, we
propose to combine the Chang-Cooper method with unconditionally positive Patankar-type time 
integration methods, since they are unconditionally positive, robust for stiff problems,
only linearly implicit, and also higher-order accurate. We describe the combined approach,
analyse it, and present a relevant numerical example demonstrating advantages compared to schemes proposed in the literature.

%% file: fokker_planck.tex
\section{Introduction}
\fpes are partial differential equations (PDEs) that are first order in time and second order in space. They can be used to describe many different social phenomena like opinion formation \cite{toscani2006kinetic,during2009boltzmann}, socio-economic phenomena \cite{furioli2017fokker}, preference formation in multi-agent societies \cite{pareschi2019hydrodynamic}, epidemic dynamics \cite{dimarco2021kinetic,zanella2023kinetic},
and spreading of fake news \cite{franceschi2022spreading}.
Furthermore, they are useful to describe processes in stochastics \cite{risken1996fokker,tome2015stochastic}, physics \cite{schuss1980singular,Elliott_2016}, biochemistry \cite{MR3609073}, and neuroscience \cite{vellmer2021fokker}.
The goal of this paper is to consider efficient numerical schemes to solve \fpes that are conservative and unconditionally positive. Those are important properties since \fpes describe a probability density which is by construction positive and conservative.  

To numerically solve the \fpeo, we follow \cite{pareschi2018structure} and use the \ccm introduced in \cite{CHANG19701} to discretise the \fpe in space. The \ccm is second-order consistent in space and designed to preserve steady states. By using the \ccmo, we obtain a semidiscretised version of the \fp equation.
To solve the resulting system of ordinary differential equations (ODEs), it is common to use the explicit Euler scheme as done in \cite{pareschi2018structure}. However, the explicit Euler scheme is not unconditionally positive and thus not unconditionally stable, see Proposition \ref{cons+pos=stab}. To avoid this, we use the \mpe and the \mprk \cite{burchard2003high,KOPECZ2018159}. These schemes are both conservative and unconditionally positive.
Such Patankar-type methods have already been used successfully for some PDEs in
\cite{ortleb2017patankar,ciallella2022arbitrary,ciallella2024high}.
Here, we apply them for the first time to \fp equations.

We apply the developed schemes to a model on opinion dynamics presented in \cite{pareschi2018structure} and derived in \cite{toscani2006kinetic}. We solve it by using different time-integration schemes and then compare the schemes with regard to computation time and numerical error.

\section{\fp Equations}
\label{sec:fpe}

We follow \cite{pareschi2018structure} and consider a \fpe of the form
\begin{equation}\label{eq:fp}
    \partial_t f(w,t) = \partial_w   \left[\mathcal{B}[f](w,t)f(w,t)+\partial_w(D(w)f(w,t)) \right]
\end{equation}
with initial condition $f(w,0)=f^0(w)$ and no-flux boundary conditions at the boundary points of $\mathcal{I}$. $t \in \mR_{\geq 0}$ denotes time, $w  \in \mathcal{I}$, $\mathcal{I}=[I_l,I_u] \subset \mR$ bounded, denotes space, $f:\mathcal{I} \times \mR_{\geq 0}\to [0, 1]$ denotes the unknown distribution function in $C^{2,1}(\mathcal{I} \times \mR_{\geq 0})$, $D:\mathcal{I} \to \mR_{\geq 0}, D\in C^2(\mathcal{I})$ denotes a diffusion function, $f^0: \mathcal{I} \to \mR_{\geq 0}$, and $f^0\in C^2(\mathcal{I})$ denotes the initial distribution. We will use the operator 
\begin{equation}
    \mathcal{B}[ f ](w, t) = \int_\mathcal{I}(w-v) f (v, t) \dd{v}
\end{equation} 
in the numerical example in Section~\ref{applications}, but other options are also possible to describe the aggregation dynamics \cite{pareschi2018structure}.
Notice that here space $w$ denotes an arbitrary quantity other than time. For example, in the applications in Section~\ref{applications}, space will be a spectrum of opinions.

The \fpe was initially used to model Brownian motion. It describes the motion of a distribution function~$f$ in the case of fluctuating macroscopic quantities. \fpes are useful when the considered variables are continuous, macroscopic, and define a small subsystem. A more detailed description of \fpes can be found in \cite{risken1996fokker,Pavliotis2014}.

To discretise problem \eqref{eq:fp} in space, we will use the \ccm which is second-order consistent in space and preserves quasi-steady states. For more information about it see \cite{pareschi2018structure,CHANG19701,ccorder}. Notice that we choose the grid points and interfaces as shown in Figure~\ref{fig:discrinspace}. $\dw$ denotes the distance between two grid points.
\begin{figure}[htbp]
\centering
    \begin{tikzpicture}
    \draw[thick] (-0.5*2.33,0) -- (5.5*2.33,0);
    \draw (2.33*0 cm,3.1pt) -- (2.33*0 cm,-3.1pt) node[anchor=north] {$w_{1}$};
    \draw (2.33*1 cm,3.1pt) -- (2.33*1 cm,-3.1pt) node[anchor=north] {$w_{2}$};
   \draw (2.33*2 cm,3.1pt) -- (2.33*2 cm,-3.1pt) node[anchor=north] {$w_{3}$};
    \draw (2.33*3 cm,3.1pt) -- (2.33*3 cm,-3.1pt) node[anchor=north] {$w_{4}$};
   \draw (2.33*4 cm,3.1pt) -- (2.33*4 cm,-3.1pt) node[anchor=north] {$w_{N-1}$};
   \draw (2.33*5 cm,3.1pt) -- (2.33*5 cm,-3.1pt) node[anchor=north] {$w_{N}$};
   \draw[black!70, thick] (-0.5*2.33 cm,3.1pt) -- (-0.5*2.33 cm,-3.1pt) node[anchor=south] {\footnotesize$w_{\frac{1}{2}}$};
   \draw[black!70, thick] (0.5*2.33 cm,3.1pt) -- (0.5*2.33 cm,-3.1pt) node[anchor=south] {\footnotesize$w_{1+\frac{1}{2}}$};
   \draw[black!70, thick] (1.5*2.33 cm,3.1pt) -- (1.5*2.33 cm,-3.1pt) node[anchor=south] {\footnotesize$w_{2+\frac{1}{2}}$};
   \draw[black!70, thick] (2.5*2.33 cm,3.1pt) -- (2.5*2.33 cm,-3.1pt) node[anchor=south] {\footnotesize$w_{3+\frac{1}{2}}$};
   \draw[black!70, thick] (3.5*2.33 cm,3.1pt) -- (3.5*2.33 cm,-3.1pt) node[anchor=south] {\footnotesize$w_{4+\frac{1}{2}}$};
   \draw[black!70, thick] (4.5*2.33 cm,3.1pt) -- (4.5*2.33 cm,-3.1pt) node[anchor=south] {\footnotesize$w_{N-\frac{1}{2}}$};
   \draw[black!70, thick] (5.5*2.33 cm,3.1pt) -- (5.5*2.33 cm,-3.1pt) node[anchor=south] {\footnotesize$w_{N+\frac{1}{2}}$};
   \draw[ultra thick] (5.5*2.33 cm,3.1pt) -- (5.5*2.33 cm,-3.1pt) node[anchor=north] {$I_u$};
   \draw[ultra thick] (-0.5*2.33 cm,3.1pt) -- (-0.5*2.33 cm,-3.1pt) node[anchor=north] {$I_l$};
\end{tikzpicture}
\caption{Discretisation in space for $N=6$ nodes.} \label{fig:discrinspace}
\end{figure}
For any continuous function $u=u(w)$, we set $u\iph:=u(w\iph)$ for $i \in \Mo:=\mathcal{M}\cup \{0\}$ with $\mathcal{M}:=\{j \in \mN| j\leq N\}$. Furthermore, $\forall t \in \mR_{\geq 0}$ we set $f_i(t):=f(w_i,t)$ for $i \in \mathcal{M}$ and $\bm{f}(t):=(f_i(t))_{i\in \mathcal{M}}$. With this, by using the \ccm we obtain the semidiscretisation
 \begin{align}\label{eq:fpfluxdi}
     \frac{\dd f_i(t)}{\dd t} = \frac{\mf\iph(t) - \mf\imh(t)}{\dw} \quad \forall t \in \mR_{\geq 0}, \quad \forall i \in \mathcal{M},
 \end{align}
 where for any $i \in \Mo$ the numerical flux $\mf\iph$ is given by
\begin{align}\label{eq:numericfluxdef}
    \mf\iph:=\mathcal{C}\iph ((1-\delta\iph)f_{i+1}+\delta\iph f_i) + D\iph \frac{f_{i+1}-f_i}{\dw},
\end{align}
with
\begin{equation}\label{eq:withoutt2}
\begin{aligned}
    \mathcal{C}\iph &= \frac{\lambda\iph D\iph}{\dw}, &
    \delta\iph &= \frac{1}{1-\exp(\lambda\iph)}+\frac{1}{\liph}, &
    \liph &= \frac{\dw(\mathcal{B}[f](w\iph)+D'\iph)}{D\iph}.
\end{aligned}
\end{equation}
For better readability, we omit writing the dependencies on $t$ in equations \eqref{eq:numericfluxdef} and \eqref{eq:withoutt2}.
As done in \cite{pareschi2018structure}, this formulation is usually solved with the \eeso. However, this is not unconditionally positive and thus a bound on the step size in time is needed. Another approach is using the implicit Euler scheme which is, however, very expensive to compute since it is in general nonlinearly implicit. Therefore, we will now propose two schemes which are unconditionally positive but only linearly implicit.

\section{Numerical Schemes for Time Integration}
To discretise \eqref{eq:fpfluxdi} in time, we consider different numerical schemes and
compare them.

First let us look at the structure--preserving properties of the numerical schemes that we consider in this paper. Therefore, from now on, if not stated otherwise, let  $n \in \mN$ and $\Delta_t \in \mR_+$ and let us set $f_i^{n}:=f_i(n\Delta_t)$ for $i \in \mathcal{M}$ and $\bm{f}^n:=(f_i^{n})_{i\in\mathcal{M}}$.

\definition A numerical scheme is
\begin{itemize}
    \item conservative, iff \begin{align}
        \sM f_i^{n+1}=\sM f_i^{n},
    \end{align}

\item unconditionally positive, iff for any positive $\bm f^n:=\bm f(t^n)$ and for any $\Delta_t\in\mR_+$ $\bm f^{n+1}$ is positive, where we call a vector positive, iff all its components are positive,
  \item stable, iff there exists $\Bar{\Delta}_t\in\mR_+$ and $\Bar{\Delta}_w\in\mR_+$ such that for any $T\in\mR_+$ there exists a constant $c_T\in\mR_+$ such that
\begin{align}
    ||\bm f^n||_{L^1} \leq c_T ||\bm f^0||_{L^1} \quad \forall 0<\dt\leq\Bar{\Delta}_t, 0<\dw\leq\Bar{\Delta}_w, 0<n\dt\leq T,
\end{align}
where $\bm f^n$ is the solution computed with the considered numerical scheme for $\dt$ and $\dw$,
\item unconditionally stable, iff it is stable for any $\Bar{\Delta}_t\in\mR_+$ and $\Bar{\Delta}_w\in\mR_+$.
\end{itemize}

Notice that there are different ways to analyse stability properties of numerical schemes for non-linear PDEs. The stability that we consider in this paper comes from the case of linear PDEs and is similar to the definition in \cite{stability}. Another way of analysing stability properties of numerical schemes for non-linear PDEs is by studying the entropy of the scheme. We do not do that in this paper, but an analysis of the entropy properties of the \ccm was done in \cite{pareschi2018structure}.  

\proposition \label{cons+pos=stab} A conservative and unconditionally positive scheme is unconditionally stable.
\proof{Let $T\in\mR_+$, $\dt,\dw \in \mR_+$ and  $0<n\dt\leq T$.
Then, \begin{align}
        ||\bm f^n||_{L^1}=\sM |f_i^n|
        =\sM f_i^n
        =\sM f_i^{n-1}=
        \dots
        =\sM f_i^1
        =\sM f^0_i      
        =\sM |f^0_i|
        =||\bm f^0||_{L^1}.        
\end{align}
\qed}

To find suitable numerical schemes in time, we follow an approach similar to \cite{burchard2003high} and consider a more general class of ODEs given in production-destruction formulation
 \begin{align}\label{eq:pd}
     \frac{\dd f_i(t)}{\dd t} = \mathfrak{P}_i(\bm f)-\mathfrak{D}_i(\bm f) \quad \forall t \in \mR_{\geq 0}, \quad \forall i \in \mathcal{M},
 \end{align}
 with initial condition $\bm f(0)=\bm f^0$     
 where 
      $t \in \mR_{\geq 0}$ denotes time,
      $\mathcal{M}=\{1,2,...,N\}$ is an index set where $N\in \mN$ denotes the number of cells at which we consider the function,
      $\bm f=\bm f(t)=(f_i(t))_{i\in \mathcal{M}}\geq 0$ denotes a vector of the unknown time-dependent functions,
      $\bm f^0=(f^0_i)_{i\in \mathcal{M}}>0$ denotes the initial vector, 
      $\mathfrak{P_i}: \mR^N \to \mR$ denotes the production rates of cell $i$, and
      $\mathfrak{D_i}: \mR^N \to \mR$ denotes the destruction rates of cell $i$.

 For $i,j \in \mathcal{M}$ let us denote the rate by which cell $j$ transforms into cell $i$ by $p_{i,j}=p_{i,j}(\bm f(t))\geq 0$ and the rate by which cell $i$ transforms into cell $j$ by $d_{i,j}=d_{i,j}(\bm f(t))\geq 0$. Clearly, $p_{i,j}=\djj$ for $j\neq i$. With this, for cell $i \in \mathcal{M}$, we can write the production terms as $\mathfrak{P}_i=\sum_{j \in \mathcal{M}}p_{i,j}$ and the destruction terms as $\mathfrak{D}_i=\sum_{j \in \mathcal{M}} d_{i,j}$.
 Since the \fpe models a distribution function which is conservative, we want to find schemes that are conservative. Therefore, we also need our system of equations to be fully conservative. Thus, we only consider the case $p_{i,i}=d_{i,i}=0$ since this leads to
 \begin{align}
     \sum_{i \in \mathcal{M}} \frac{\dd f_i}{\dd t}= \sum_{i \in \mathcal{M}} (\mathfrak{P}_i-\mathfrak{D}_i) = \sum_{i \in \mathcal{M}} \sum_{j \in \mathcal{M}} (\pii-\dii) =  \sum_{i \in \mathcal{M}} (p_{i,i}-d_{i,i}) =0, 
 \end{align}
 which shows that the sum over the component functions is constant in time, and thus the system of equations is fully conservative. 

Notice that we can write the semi-discretised formulation \eqref{eq:fpfluxdi} of our initial \fpe \eqref{eq:fp} in the formulation of problem \eqref{eq:pd}. One way of doing this is by setting for $i\in \mathcal{M}\setminus \{1,N\}$
 \begin{align}\label{eq:choice_p_d}
     \begin{split}
            &p_{i,i+1}(\bm f)=d_{i+1,i}:=\frac{\max\left(0,\Tilde{\mathcal{C}}_{i+1/2}\right)\left((1-\delta\iph)f_{i+1}+\delta\iph f_i\right)
                         +\frac{D_{i+1/2}f_{i+1}}{\Delta_w}}{\Delta_w}\geq 0,\\
            &p_{i,i-1}(\bm f)=d_{i-1,i}:=\frac{-\min\left(0,\Tilde{\mathcal{C}}_{i-1/2}\right)\left((1-\delta\imh)f_{i}+\delta\imh f_{i-1}\right)
            +\frac{D_{i-1/2}f_{i-1}}{\Delta_w}}{\Delta_w}\geq 0,\\
            &p_{i,j}(\bm f)=d_{j,i}(\bm f):=0 \quad \forall \quad j\in \mathcal{M}\setminus\{i-1,i+1\}.
     \end{split}
 \end{align}

We consider the unconditionally positive \textbf{modified Patankar-Euler scheme (MPE)}
 \begin{align}\label{eq:mpe}
     f_i^{n+1} = f_i^{n} + \Delta_t\left(\sMj \pii(\bm f^n)\frac{f_j^{n+1}}{f_j^{n}}-\sMj \dii(\bm f^n)\frac{f_i^{n+1}}{f_i^{n}}\right) \quad \forall i \in \mathcal{M},
 \end{align} and the \textbf{modified Patankar-Runge-Kutta scheme (MPRK)}
 \begin{align}\label{eq:mprk}
\begin{split}
         \Tilde{f_i} &= f_i^{n} + \Delta_t\left(\sMj \pii(\bm f^n)\frac{\Tilde{f_j}}{f_j^{n}}-\sMj \dii(\bm f^n)\frac{\Tilde{f_i}}{f_i^{n}}\right) \quad \forall i \in \mathcal{M},\\
        f_i^{n+1} &= f_i^{n} + \frac{\Delta_t}{2}\left(\sMj \left(\pii(\bm f^n)+\pii(\Tilde{\bm f})\right)\frac{f_j^{n+1}}{\Tilde{f_j}}-\sMj \left(\dii(\bm f^n)+\dii(\Tilde{\bm f})\right)\frac{f_i^{n+1}}{\Tilde{f_i}}\right) \quad \forall i \in \mathcal{M}.
\end{split}
 \end{align}
As shown in \cite{burchard2003high,KOPECZ2018159}, these schemes are unconditionally positive and conservative. The \mpe is first-order accurate and the \mprk is second-order accurate, see \cite{KOPECZ2018159}.
Notice that since the \mper is unconditionally positive, $\Tilde{\bm f}$ is positive for a positive $\bm f^n$, and thus the second step of the \mprkr is feasible.
We do not consider other members of the family of second-order Patankar-type methods
of \cite{KOPECZ2018159} since they can have reduced performance for solutions near
zero \cite{torlo2022issues}.

Since the \mper and the \mprkr are unconditionally positive while the explicit Euler scheme is only positive and stable under a restriction on the time step we can take larger time step sizes when solving a problem with the \mpe or the \mprko. However, when using the \mpe we have to solve one system of linear equations and when using the \mprk two systems of linear equations in each time step, which means that a time step of one of those schemes requires more computing power than a time step of the explicit Euler scheme or the Heun scheme, respectively.

These Patankar-type methods are not only unconditionally positive and
conservative but also work well for stiff problems. This has been observed
numerically in \cite{KOPECZ2018159}. 
This property is also important since the semidiscretisations of the 
\fpe yield a parabolic time step size restriction for explicit methods. A new approach to analyse the stability
of these methods has been proposed recently \cite{izgin2022lyapunov}.

The \ccm is designed to preserve steady states. Standard time integration
methods such as the \ees typically preserve steady states, too, but more
involved methods like IMEX schemes do not necessarily share this property
\cite{chertock2015steady}. However, the Patankar-type methods are 
well-suited for these applications since 
they preserve steady states \cite{torlo2022issues}.

\section{Application}\label{applications}
In this section, we use the schemes introduced before to solve one model that can be described by \fpeso. The model presented here is from \cite{pareschi2018structure}. For a derivation and more thorough analysis of the model, see \cite{toscani2006kinetic}.
In this model, $f(w,t)$ describes the ratio of people having opinion $w$ at time $t\in\mR_{\geq 0}$, where the spectrum of opinions is given by the interval $\mathcal{I}=[-1,1]$. The \fpe is determined by 
\begin{align}\label{eq:B}
\begin{split}
       \mathcal{B}[ f ](w, t) =\int_\mathcal{I}(w-v) f (v, t) \dd{v}, \qquad
       D(w)=\frac{\sigma^2}{2}(1-w^2)^2,
\end{split}
\end{align}
for a parameter $\sigma\in \mR$. In the following numerical tests, we choose $\sigma=\sqrt{0.2}$.
Like in our initial problem \eqref{eq:fp}, we consider the problem for no-flux boundary conditions. As initial distribution we choose
\begin{align}
    f^0(w)=\beta\left(\exp\left(-30\left(w+\frac{1}{2}\right)^2\right)+\exp\left(-30\left(w-\frac{1}{2}\right)^2\right)\right),
\end{align}
where the constant $\beta\in \mR_+$ is chosen such that $\int_\mathcal{I}f^0\dd w=1$. 

As explained in \cite{toscani2006kinetic}, under those conditions, it is possible to analytically obtain the stationary solution $f_\infty(w) = \lim\limits_{t\to\infty}f(w,t)$
\begin{align}
    f_\infty(w)=\frac{K}{(1-w^2)^2}\left(\frac{1+w}{1-w}\right)^\frac{u}{2\sigma^2}\exp\left(-\frac{1-uw}{\sigma^2(1-w^2)}\right),
    \quad 
    u := \int_\mathcal{I}w f^0(w) \dd w.
\end{align}

We solve the \fpe by using the \ccm to obtain the numerical fluxes, and then use the \mpeo, the \mprk to solve the problem in its semidiscrete formulation. We compare the Patankar methods with the \ees and the Heun scheme since the Patankar schemes are based on these methods. As a fully implicit method, we also use the implicit Euler method.
We consider the solution over the time interval $[0,10]$ and, if not stated otherwise, we use $N=80$ grid points for the discretisation in space. Since we consider the problem for the space interval $\mathcal{I}=[-1,1]$, for $N=80$ we get $\dw=0.025$. For the Patankar schemes we use the Julia \cite{bezanson2015julia} package PositiveIntegrators.jl \cite{kopecz2023positiveintegrators}. OrdinaryDiffEq.jl \cite{rackauckas2017differentialequations} is applied for the classical Runge-Kutta methods. The source code used for this paper can be found in our reproducibility repository \cite{bartel2024structureRepro}.

Heun's second-order method is strong stability preserving (SSP), i.e., it preserves convex stability properties of the explicit Euler method such as positivity \cite{gottlieb2011strong}.
Concretely, if the explicit Euler method is positive under the time step restriction $\Delta_t \le \Delta_{t, \mathrm{FE}}$, an SSP method is positive for $\Delta_t \le \mathcal{C} \Delta_{t, \mathrm{FE}}$, where
$\mathcal{C}$ is the SSP coefficient of the method. The implicit Euler method is unconditionally positive since it is SSP with $\mathcal{C} = \infty$ \cite[Chapter~7]{gottlieb2011strong}. Heun's method has the SSP coefficient $\mathcal{C} = 1$ \cite[Section~2.4.1]{gottlieb2011strong}, i.e., it is positive under the same time step size restriction as the explicit Euler method.

\begin{figure}[htbp]
\begin{minipage}{72mm}
\centering
{\includegraphics[width=1\linewidth]{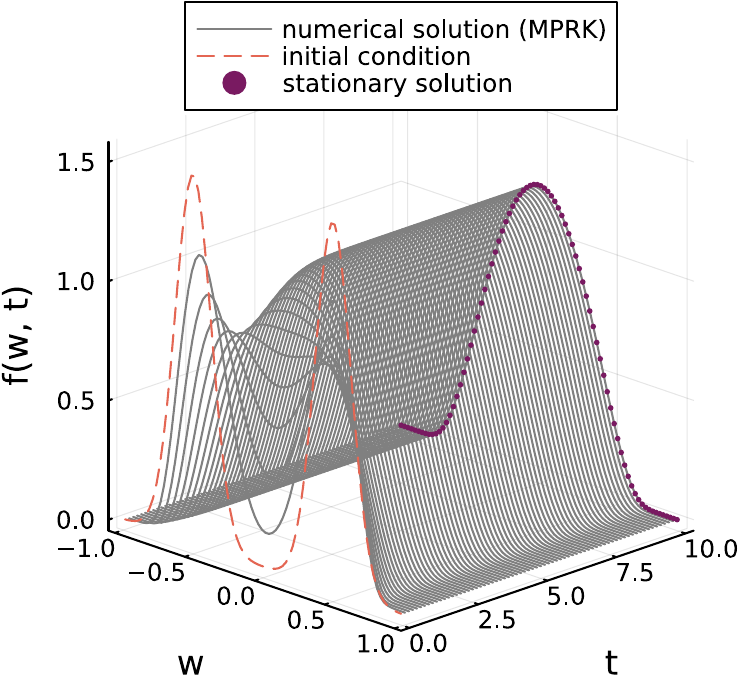}}
  \caption{Numerical solution computed with \ccm and \mprk for \tww.}
\label{ccmmprk1}
\end{minipage}
\hfil
\begin{minipage}{72mm}
\centering
{\includegraphics[width=1\linewidth]{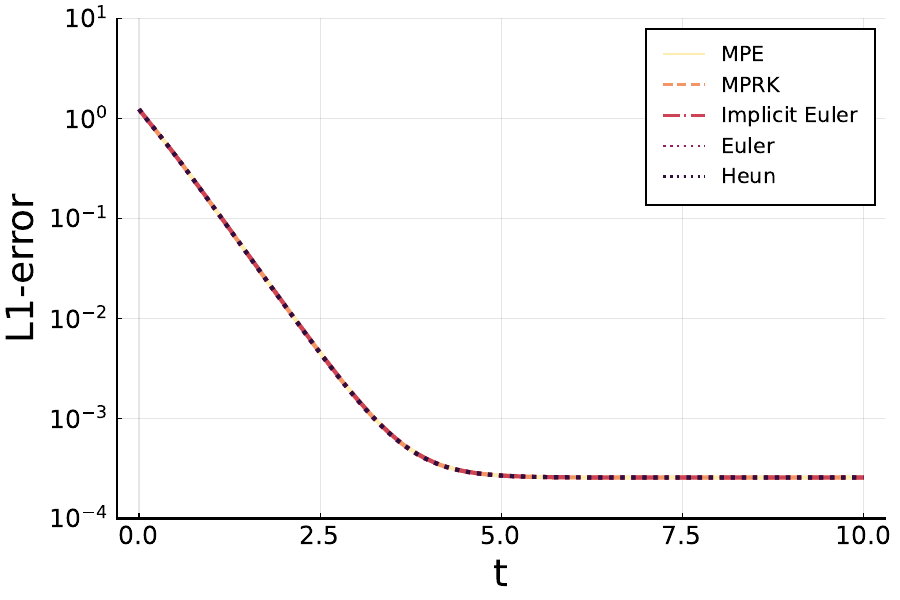}}
\caption{$L^1$-error compared to stationary solution for different schemes and \tww.}
  \label{L1-errordtdw^21}
\end{minipage}
\end{figure}

\begin{figure}[htbp]
\begin{minipage}{72mm}
\centering
{\includegraphics[width=1\linewidth]{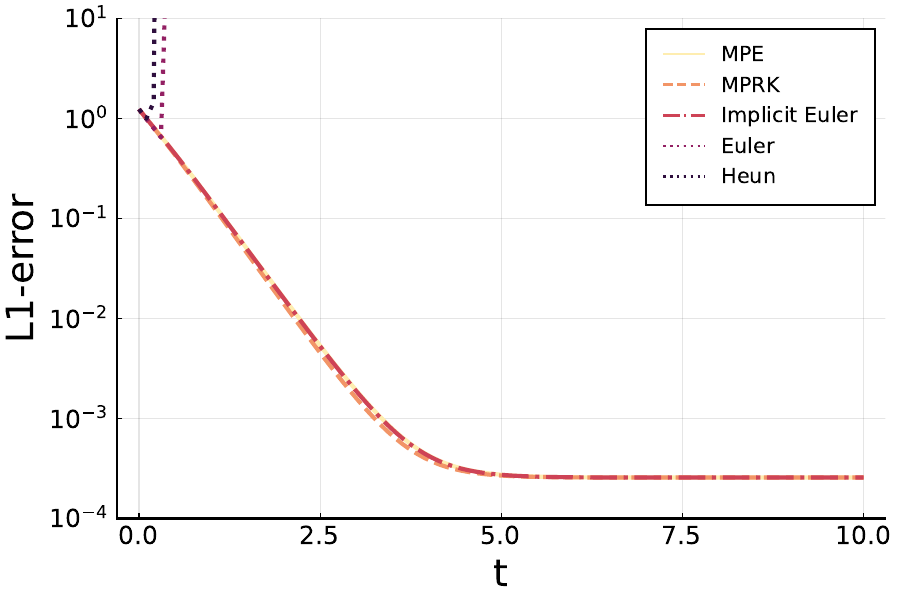}}
      \caption{$L^1$-error compared to stationary solution for different schemes and $\dt=\dw$.}
  \label{L1-errordtdw1}
\end{minipage}
\hfil
\begin{minipage}{72mm}
\centering
{\includegraphics[width=1\linewidth]{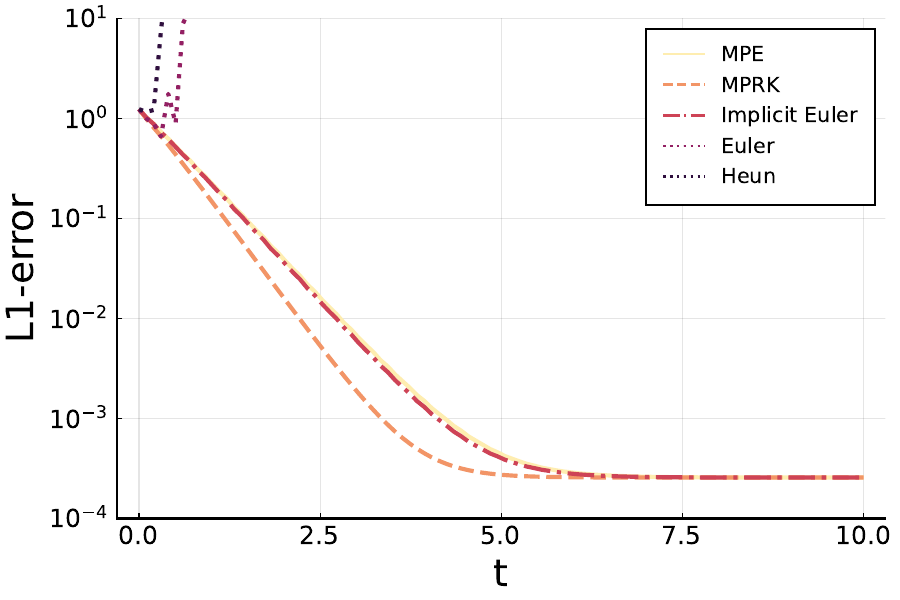}}
 \caption{$L^1$-error compared to stationary solution for different schemes and $\dt=10\dw$.}
  \label{L1-errordt10dw1}
\end{minipage}
\end{figure}

In Figure~\ref{ccmmprk1}, we can see the numerical solution computed with the \ccm combined with the \mprk for \tww. We can see that it converges to the stationary solution. The behaviour of the solution computed for this choice of $\dt$ computed with the \ccm combined with the \mpeo, the \mprko, the implicit Euler scheme, the \eeso, and the Heun scheme look the same which is why we here only display one. The others can be computed with our code \cite{bartel2024structureRepro}. That the five schemes lead to a similar numerical solution can also be seen in Figure~\ref{L1-errordtdw^21} and occurs since for that choice of $\dt$ the positivity conditions for the explicit Euler and the Heun scheme are satisfied.
Figures \ref{L1-errordtdw^21}, \ref{L1-errordtdw1} and \ref{L1-errordt10dw1} display the $L^1$-error of the numerical solution compared to the stationary solution, where the numerical solution is computed for the five considered schemes and for three different choices of $\dt$. In Figure~\ref{L1-errordtdw^21}, where we choose \tww like in Figure~\ref{ccmmprk1}, we can see that the $L^1$-error decreases similarly for all five schemes and converges to around $0.00077$. In Figures \ref{L1-errordtdw1} and \ref{L1-errordt10dw1}, however, the choices $\dt=\dw$ and $\dt=10\dw$ do not satisfy the positivity conditions for the explicit Euler and Heun scheme. Thus, the $L^1$-error for these schemes goes to infinity, which shows that they are not stable for those choices of $\dt$, and due to Proposition \ref{cons+pos=stab} we know that that is due to the lack of positivity. Since the implicit Euler scheme, the \mpe and the \mprk are unconditionally positive, for these three schemes the $L^1$-error converges to a fixed value that is depending on the discretisation in space. Moreover, Figure~\ref{ccmmprk1} shows that for large choices of $\dt$, the implicit Euler scheme and the \mpe converge slower to the stationary solution than for the other $\dt$ while the \mprk converges similarly for all the considered choices of $\dt$.

\begin{figure}[htbp]
\begin{minipage}{72mm}
\centering
{\includegraphics[width=1\linewidth]{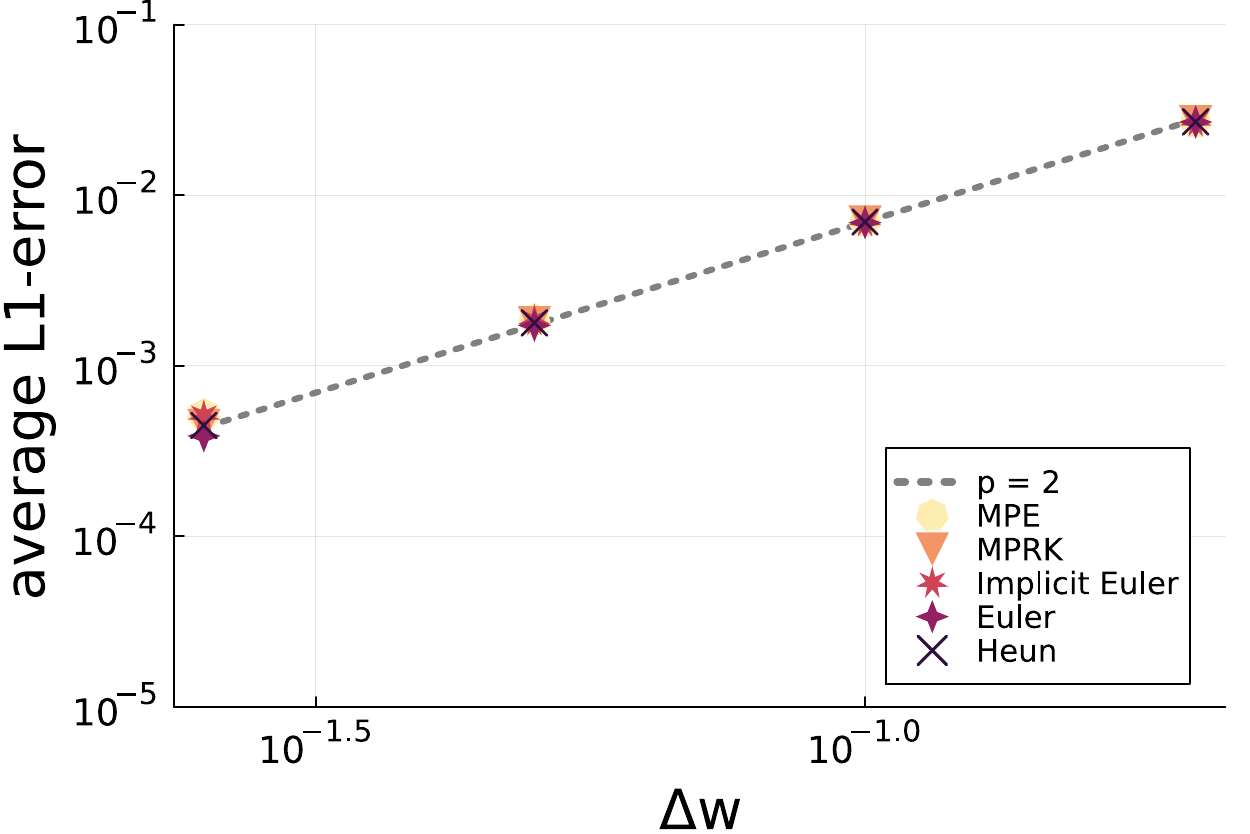}}
      \caption{Experimental order of convergence in space.}
  \label{oc_space}
\end{minipage}
\hfil
\begin{minipage}{72mm}
\centering
{\includegraphics[width=1\linewidth]{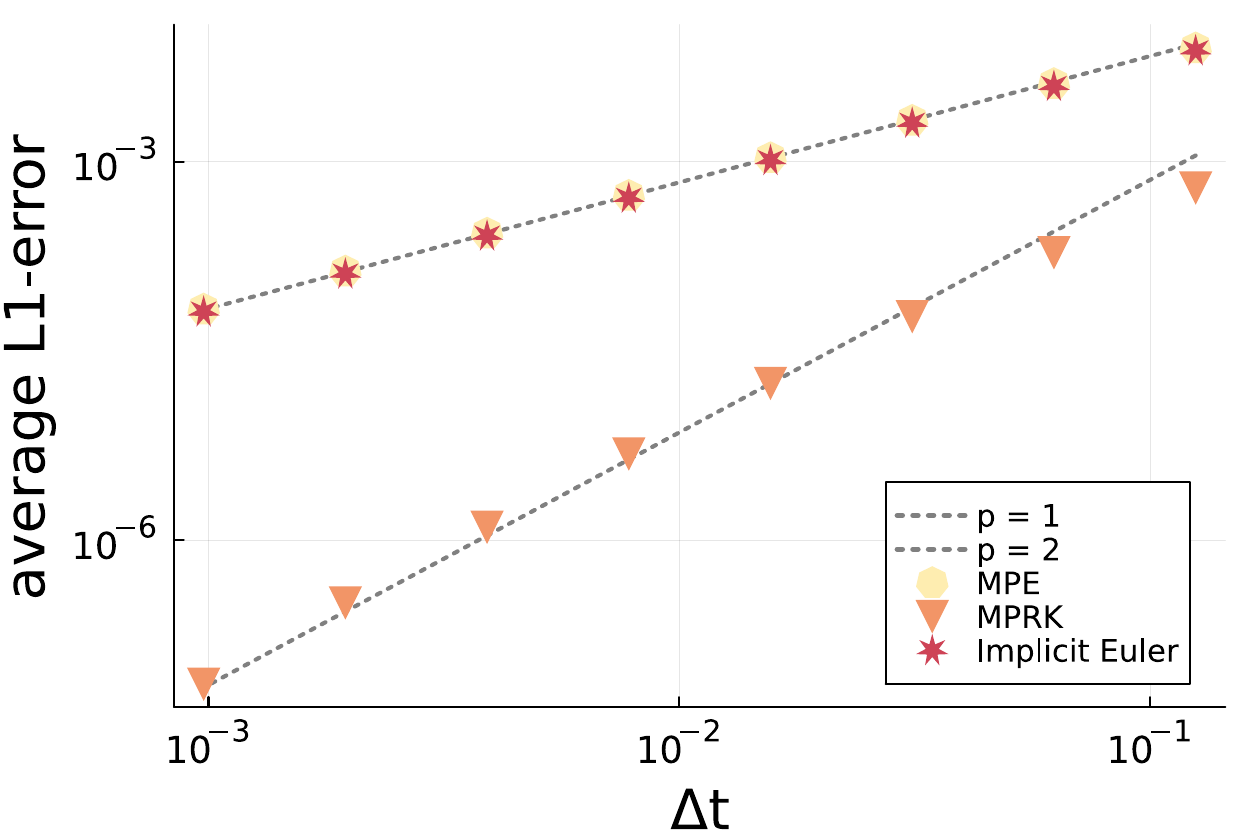}}
 \caption{Experimental order of convergence in time.}
  \label{oc_time}
\end{minipage}
\end{figure}

Now let us look at the experimental order of convergence. We compute the $L^1$-error averaged over time with respect to a reference solution, which is computed with $N=640$ grid points and the \ccm combined with the explicit Euler scheme for \twwo. Due to our choice of the grid (see Figure~\ref{fig:discrinspace}), this does not give us the values of the reference solution at the grid points where we know the values of the numerical solution. Therefore, we use cubic spline interpolation from the Julia package Interpolations.jl \cite{kittisopikul2022interpolations} to interpolate the reference solution with a twice continuously differentiable piece-wise cubic polynomial. Figure~\ref{oc_space} displays the experimental order of convergence in space, and we can see that for all the five different schemes we used for the time steps, our scheme is of second order in space. Figure~\ref{oc_time} shows the experimental order of convergence in time for the \mpeo, the \mprk and the implicit Euler scheme. We can see that the \mpe and the implicit Euler scheme are first-order convergent and the \mprk reaches an order of convergence of two. To avoid errors due to interpolation, we use $N = 160$ points in space for both the reference solution and the numerical solutions. The reference solution uses the Heun method with \twwo.

Next, we study the computation times of the different methods. Therefore, we compute the solution with the five schemes for $N = 80$ grid points and different values of $\dt$. For the computations, we use a laptop with an Intel\textregistered\ Core i5-1345U processor with 12 cores and a RAM of 16~GB.
First, we compare the total CPU time needed to compute the solution. Hence, we measure the time needed to evaluate the code to compute the solution over the time interval $[0,10]$ for the different schemes. 
Table~\ref{comp_time} presents the time needed for one evaluation of the code averaged over five runs for each scheme. As expected, the computing times decrease for increasing $\dt$. Moreover, Table~\ref{comp_time} shows that for all considered choices of $\dt$, the \mpe and \mprk take roughly 50~\% longer than their corresponding classical Runge-Kutta methods. Due to the need to solve nonlinear systems in each time step, the implicit Euler methods takes much more time than the other schemes.
However, due to the condition on the time step for the \eeso, we can use the \ees only for small time steps. On the other hand, by choosing a large time step, the computation time reduces significantly. Thus, the fastest choice is the \mpe with a large time step, which is just slightly faster than the \mprk when considering large time steps.

\begin{table}
    \caption{Mean total CPU time and standard deviation when computing $\frac{10}{\dt}$ time steps (averaged over five runs).}
    \label{comp_time}
    \footnotesize
\begin{tabular}{cccccc}
  \toprule
  \textbf{$\Delta_t$} & \textbf{MPE} & \textbf{MPRK} & \textbf{Implicit Euler} & \textbf{Euler} & \textbf{Heun} \\\midrule
  $\frac{\Delta_w^{2.5}}{2\sigma^2}$ & (239.59 $\pm$ 7.57)~ms & (485.12 $\pm$ 5.33)~ms & 14.33~s $\pm$ 61.21~ms & (166.28 $\pm$ 0.63)~ms & (325.72 $\pm$ 1.23)~ms \\
  $\frac{\Delta_w^2}{2\sigma^2}$ & (37.77 $\pm$ 0.05)~ms & (77.1 $\pm$ 0.21)~ms & 2.26~s $\pm$ 3.98~ms & (25.94 $\pm$ 0.05)~ms & (51.43 $\pm$ 0.14)~ms \\
  $\Delta_w$ & (2.34 $\pm$ 0.02)~ms & (4.81 $\pm$ 0.04)~ms & (142.02 $\pm$ 0.74)~ms & -- & -- \\
  $\frac{\Delta_w}{2\sigma^2}$ & (0.92 $\pm$ 0.02)~ms & (1.9 $\pm$ 0.02)~ms & (56.42 $\pm$ 0.16)~ms & -- & -- \\
  $10\Delta_w$ & (0.22 $\pm$ 0.003)~ms & (0.47 $\pm$ 0.03)~ms & (14.18 $\pm$ 0.03)~ms & -- & -- \\\bottomrule
\end{tabular}
\end{table}

\begin{figure}[tbp]
\begin{minipage}{72mm}
\centering
   {\includegraphics[width=1.15\linewidth]{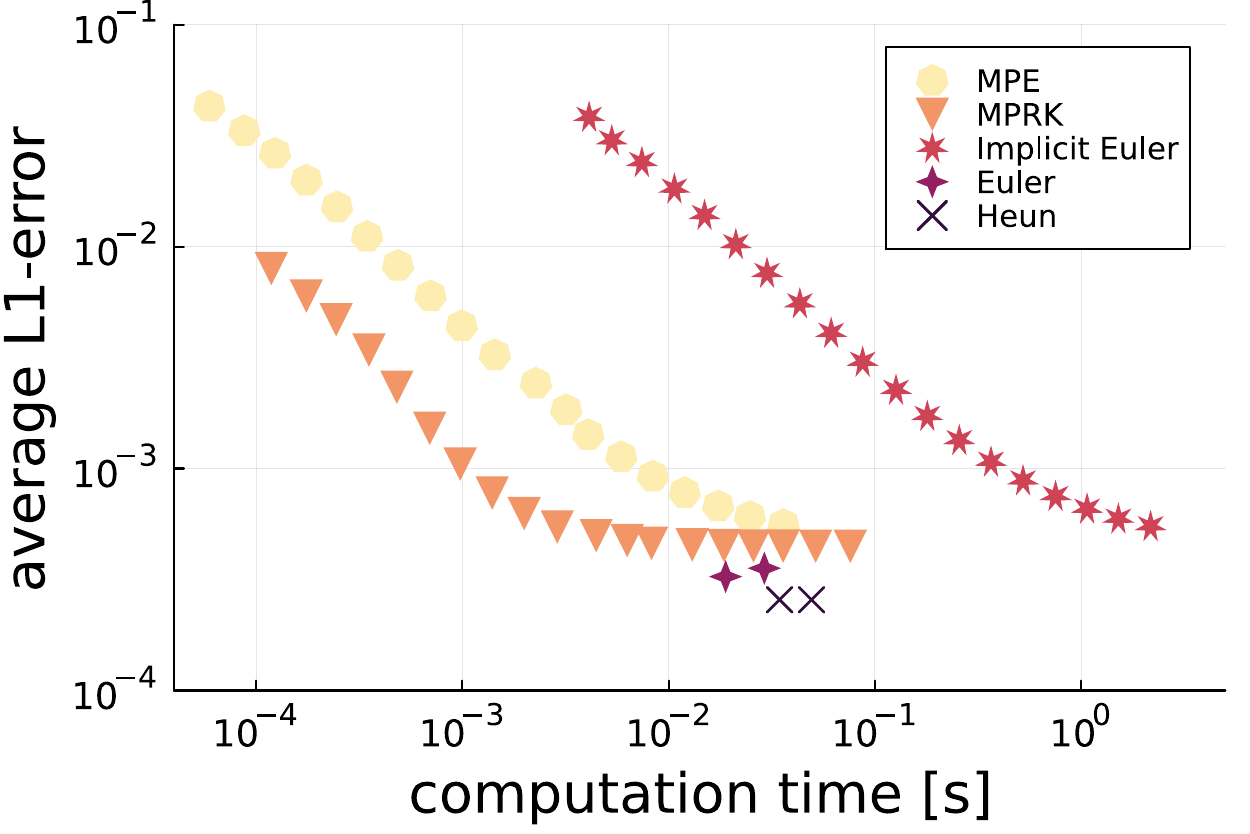}}
  \caption{Computation time vs. average $L^1$-error (median over five runs).}
  \label{compvsl1av}
\end{minipage}
\hfil
\begin{minipage}{72mm}
\centering
 {\includegraphics[width=1\linewidth]{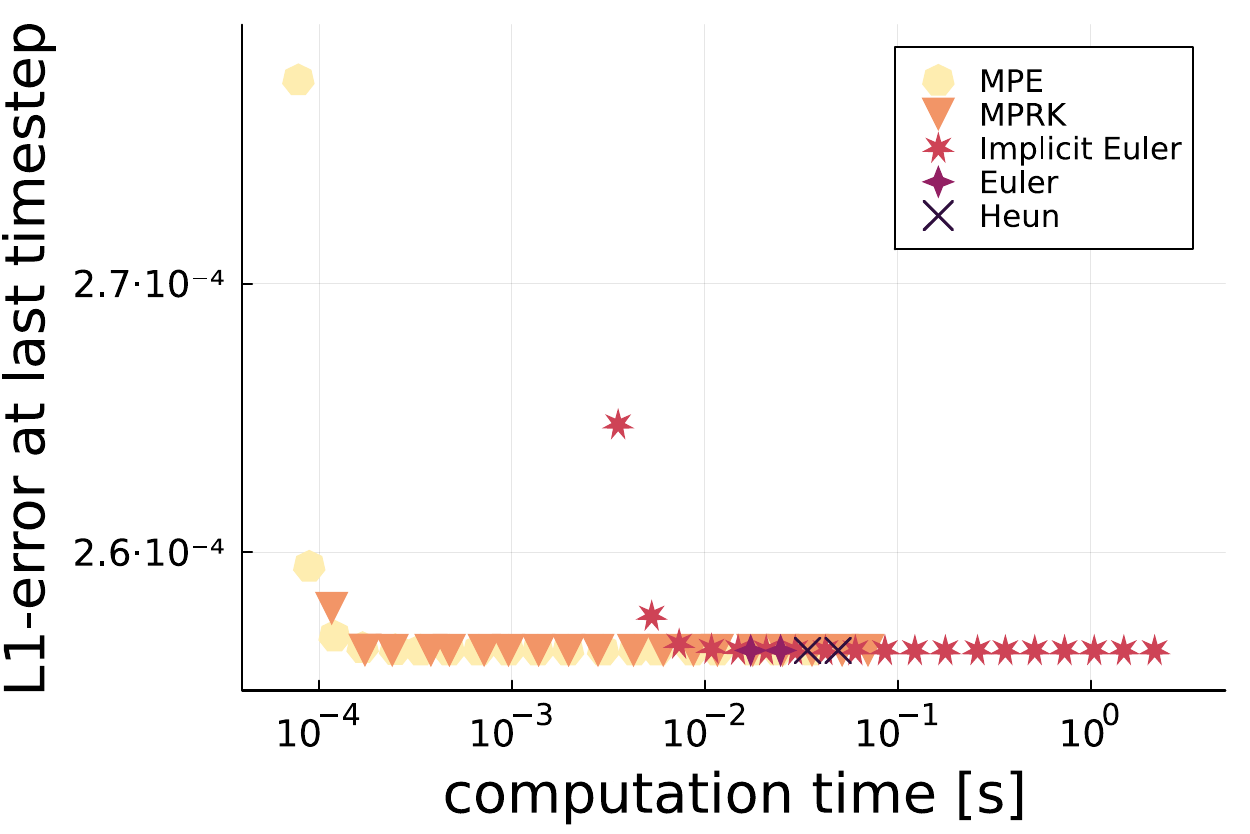}}
  \caption{Computation time vs.\ $L^1$-error at the final time (median over five runs).}
  \label{compvsl1_stat}
\end{minipage}
\end{figure}

Secondly, we want to compare the relation between the computation time needed to compute the solution and the $L^1$-error of the solution for the different schemes. Therefore, we measure the wall time needed to compute the solution for $\dt\in \mathcal{T}:=\{0.7^k|k\in\mN_0, k< 19\}$ with the \ccm combined with the \mpeo, the \mprko, the implicit Euler scheme, the \eeso, and the Heun scheme. 
Figure~\ref{compvsl1av} displays the computation time needed to compute the solution for $\dt\in\mathcal{T}$ versus the average $L^1$-error of the numerical solution with the reference solution.

Figure~\ref{compvsl1_stat} shows the computation time needed to compute the solution for $\dt\in\mathcal{T}$ versus the $L^1$-error of the numerical solution and the stationary solution at the last time step, i.e., how well the numerical solution converges to the stationary solution in the considered time interval. In both figures, we see that the $L^1$-error for the solution computed with the \ees and the Heun scheme is only displayed twice. This is because the \ees and the Heun scheme are only stable for small $\dt$. We can see in Figure~\ref{compvsl1av} that with increasing computation time, the average $L^1$-error decreases. Furthermore, Figure~\ref{compvsl1av} shows that, especially for short computation times, the average $L^1$-error of the solution computed with the \mprk is less than for the solution computed with the \mpeo. Notably, the \mprk can reach almost the same $L^1$-error as the explicit Euler and Heun scheme in less time. In Figure~\ref{compvsl1_stat} we can observe that for $\dt\leq 0.7^{3}$ the $L^1$-error of the numerical solution at the last time step and the stationary solution is approximately the same for the solutions computed with the \mpe and the \mprko. For $\dt= 0.7^{17}$ and $\dt= 0.7^{18}$ the $L^1$-error at the last time step of the solution obtained by the \ees is also similar. For $\dt= 0.7^0$ the $L^1$-error at the last time step of the solution obtained by the \mprk is higher than for the other choices of $\dt$. It is, however, smaller than the $L^1$-error at the last time step of the solution obtained by the \mpe for $\dt= 0.7^0$.
The implicit Euler method behaves very similar to the \mpeo, but with a computational cost that is more than an order of magnitude bigger.

From Table~\ref{comp_time}, Figure~\ref{compvsl1av}, and Figure~\ref{compvsl1_stat} we can conclude that it is more efficient to use the \mpe or \mprk instead of the \ees since, even though the \ees is the fastest per time step, the \ees is not stable for choices of $\dt$ where the other schemes are faster and have a small $L^1$-error. However, their average $L^1$-error is bigger than the $L^1$-error of the \ees for all the considered choices of $\dt$. 
For very large values of $\dt$ it is better to use the \mprk since the computation time needed is only a little more than when using the \mpeo, but the computed solution approximates the reference solution better and converges closer to the stationary solution. Figure~\ref{compvsl1_stat} shows that the Patankar methods are especially well-suited if one is only interested in the stationary solution since large time steps already suffice for a good approximation. Otherwise, the \mpe is a good choice since its computation time is less than for the \mprk and for rather small values of $\dt$ its $L^1$-error is similar to the $L^1$-error of the solution computed with the \mprko.

\section{Conclusion}
We discussed the possibility of solving \fpes with the \ccm combined with either the \mpe or the \mprko. Those schemes are conservative, unconditionally positive, and preserve steady states.
These are essential properties when solving \fpes since they describe a probability density and, thus, are positive and conservative. We also solved a specific \fpe modelling opinion dynamics with the \ccm combined with the \mpeo, the \mprko, and the explicit and implicit Euler scheme and the Heun scheme. We compared the performance of the different schemes, which showed the advantages of the Patankar-type schemes.

Our results suggest that the unconditional positivity of the Patankar-type methods comes at the price of reduced accuracy compared to their baseline explicit Runge-Kutta methods.
Thus, it is more efficient to use the standard methods when the time step size is so small that positivity can be ensured.
However, the Patankar-type methods work well for stiff problems and large time step sizes.
In this regime, they can be significantly more efficient than classical (explicit or fully implicit) Runge-Kutta methods.

%% file: acknowledgments.tex
HB and JL acknowledge the support by the Deutsche Forschungsgemeinschaft (DFG, 
German Research Foundation) within the Research Training Group GRK 2583 
``Modeling, Simulation and Optimisation of Fluid Dynamic Applications''.
HR was supported by the Deutsche Forschungsgemeinschaft (DFG,
German Research Foundation, project number 513301895) and the 
Daimler und Benz Stiftung (Daimler and Benz foundation, project number 32-10/22).

%% file: references.bib
@misc{stability,
  author        = {Koev, Plamen},
  title         = {Numerical Methods for Partial Differential Equations},
  year          = {2012},
url={https://math.mit.edu/~plamen/18.336/math336.pdf}
}

@Inbook{Pavliotis2014,
author="Pavliotis, Grigorios A.",
title="The Fokker--Planck Equation",
bookTitle="Stochastic Processes and Applications: Diffusion Processes, the Fokker--Planck and Langevin Equations",
year="2014",
publisher="Springer New York",
address="New York, NY",
pages="87--137",
doi="https://doi.org/10.1007/978-1-4939-1323-7_4"
}

@article{franceschi2022spreading,
  title={Spreading of fake news, competence and learning: kinetic modelling and numerical approximation},
  author={Franceschi, Jonathan and Pareschi, Lorenzo},
journal = {Philosophical Transactions of the Royal Society A: Mathematical, Physical and Engineering Sciences},
  volume={380},
  number={2224},
  year={2022},
  publisher={The Royal Society},
doi={https://doi.org/10.1098/rsta.2021.0159}
}

@article{during2009boltzmann,
  title={Boltzmann and Fokker--Planck equations modelling opinion formation in the presence of strong leaders},
  author={D{\"u}ring, Bertram and Markowich, Peter and Pietschmann, Jan--Frederik and Wolfram, Marie--Therese},
  journal={Proceedings of the Royal Society A: Mathematical, Physical and Engineering Sciences},
  volume={465},
  number={2112},
  pages={3687--3708},
  year={2009},
  publisher={The Royal Society Publishing},
doi={https://doi.org/10.1098/rspa.2009.0239}
}

@article{furioli2017fokker,
  title={Fokker--Planck equations in the modeling of socio--economic phenomena},
  author={Furioli, Giulia and Pulvirenti, Ada and Terraneo, Elide and Toscani, Giuseppe},
  journal={Mathematical Models and Methods in Applied Sciences},
  volume={27},
  number={01},
  pages={115--158},
  year={2017},
  publisher={World Scientific},
doi={https://doi.org/10.1142/S0218202517400048}
}

@article{pareschi2019hydrodynamic,
  title={Hydrodynamic models of preference formation in multi--agent societies},
  author={Pareschi, Lorenzo and Toscani, Giuseppe and Tosin, Andrea and Zanella, Mattia},
  journal={Journal of Nonlinear Science},
  volume={29},
  pages={2761--2796},
  year={2019},
  publisher={Springer},
doi={https://doi.org/10.1007/s00332-019-09558-z}
}

@article{toscani2006kinetic,
  title={Kinetic models of opinion formation},
  author={Toscani, Giuseppe},
  journal={Communications in Mathematical Sciences},
  volume={4},
  number={3},
  pages={481--496},
  year={2006},
  publisher={International Press of Boston},
doi = {10.4310/CMS.2006.v4.n3.a1}
}

@article{ccorder,
  title={Analysis of the Chang--Cooper discretization scheme for a class of Fokker--Planck equations},
  author={Masoumeh Mohammadi and Alfio Borz{\`i}},
  journal={Journal of Numerical Mathematics},
  year = {2015},
number = {3},
  volume={23},
  pages={271--288},
    doi = {doi:10.1515/jnma-2015-0018}
}

@article{burchard2003high,
title = {A high--order conservative Patankar--type discretisation for stiff systems of production--destruction equations},
journal = {Applied Numerical Mathematics},
volume = {47},
number = {1},
pages = {1--30},
year = {2003},
doi = {https://doi.org/10.1016/S0168-9274(03)00101-6},
author = {Hans Burchard and Eric Deleersnijder and Andreas Meister},
}

@article{pareschi2018structure,
  title={Structure preserving schemes for nonlinear Fokker--Planck equations and applications},
  author={Pareschi, Lorenzo and Zanella, Mattia},
  journal={Journal of Scientific Computing},
  volume={74},
  pages={1575--1600},
  year={2018},
  publisher={Springer},
  doi={https://doi.org/10.1007/s10915-017-0510-z}
}

@article{KOPECZ2018159,
title = {On order conditions for modified Patankar--Runge--Kutta schemes},
journal = {Applied Numerical Mathematics},
volume = {123},
pages = {159--179},
year = {2018},
doi = {https://doi.org/10.1016/j.apnum.2017.09.004},
author = {Stefan Kopecz and Andreas Meister}
}

@article{CHANG19701,
title = {A practical difference scheme for Fokker--Planck equations},
journal = {Journal of Computational Physics},
volume = {6},
number = {1},
pages = {1-16},
year = {1970},
doi = {https://doi.org/10.1016/0021-9991(70)90001-X},
author = {Julius S. Chang and Gilbert E. Cooper},
}

@article{dimarco2021kinetic,
  title={Kinetic models for epidemic dynamics with social heterogeneity},
  author={Dimarco, Giacomo
and Perthame, Benoit
and Toscani, Giuseppe
and Zanella, Mattia},
  journal={Journal of Mathematical Biology},
  volume={83},
  number={4},
  year={2021},
  publisher={Springer},
doi={https://doi.org/10.1007/s00285-021-01630-1}

}

@article{zanella2023kinetic,
  title={Kinetic models for epidemic dynamics in the presence of opinion polarization},
  author={Zanella, Mattia},
  journal={Bulletin of Mathematical Biology},
  volume={85},
  number={36},
  year={2023},
  publisher={Springer},
doi={10.1007/s11538-023-01147-2},
}

@book{risken1996fokker
,	title	= {The Fokker--Planck equation}
, subtitle  = {Methods of solution and applications}
,	author	= {Risken, Hannes}
,   series  = {Springer Series in Synergetics}
,	volume	= {18}
,	year	= {1996}
,	publisher	= {Springer}
,  address = {Berlin, Heidelberg},
DOI={https://doi.org/10.1007/978-3-642-61544-3},
}

@book{tome2015stochastic,
  title={Stochastic dynamics and irreversibility},
  author={Tom{\'e}, T{\^a}nia and de Oliveira, M{\'a}rio J.},
  year={2015},
  publisher={Springer},
doi = {https://doi.org/10.1007/978-3-319-11770-6},
}

@article{schuss1980singular,
author = {Schuss, Zeev},
title = {Singular Perturbation Methods in Stochastic Differential Equations of Mathematical Physics},
journal = {SIAM Review},
volume = {22},
number = {2},
pages = {119--155},
year = {1980},
doi = {https://doi.org/10.1137/1022024},
}

@article {MR3609073,
    AUTHOR = {Schnoerr, David and Sanguinetti, Guido and Grima, Ramon},
     TITLE = {Approximation and inference methods for stochastic biochemical
              kinetics---a tutorial review},
   JOURNAL = {J. Phys. A},
  FJOURNAL = {Journal of Physics. A. Mathematical and Theoretical},
    VOLUME = {50},
      YEAR = {2017},
    NUMBER = {9},
     PAGES = {093001, 60},
 MRCLASS = {60J28 (60J70 62F15 92C45)},
  MRNUMBER = {3609073},
       DOI = {https://doi.org/10.1088/1751-8121/aa54d9},
}

@article{Elliott_2016,
	doi = {10.1103/physreva.94.043840}, 
	url = {https://doi.org/10.1103%2Fphysreva.94.043840},
	year = 2016,
	month = {oct},
	publisher = {American Physical Society ({APS})},
	volume = {94},
	number = {4},
	author = {Matthew Elliott and Eran Ginossar},
	title = {Applications of the Fokker--Planck equation in circuit quantum electrodynamics},
	journal = {Physical Review A}
}

@article{vellmer2021fokker,
  title={Fokker--Planck approach to neural networks and to decision problems},
  author={Vellmer, Sebastian and Lindner, Benjamin},
  journal={The European Physical Journal Special Topics},
  volume={230},
  number={14},
  pages={2929--2949},
  year={2021},
  publisher={Springer},
  doi={https://doi.org/10.1140/epjs/s11734-021-00172-3}
}

@article{torlo2022issues,
  title={Issues with Positivity-Preserving {P}atankar-type Schemes},
  author={Torlo, Davide and \"Offner, Philipp and Ranocha, Hendrik},
  journal={Applied Numerical Mathematics},
  year={2022},
  month={08},
  volume={182},
  pages={117--147},
  doi={10.1016/j.apnum.2022.07.014}
}

@article{izgin2022lyapunov,
  title={On {L}yapunov stability of positive and conservative time integrators and application to second order modified {P}atankar--{R}unge--{K}utta schemes},
  author={Izgin, Thomas and Kopecz, Stefan and Meister, Andreas},
  journal={ESAIM: Mathematical Modelling and Numerical Analysis},
  volume={56},
  number={3},
  pages={1053--1080},
  year={2022},
  publisher={EDP Sciences},
  doi={10.1051/m2an/2022031}
}

@article{chertock2015steady,
  title={Steady state and sign preserving semi-implicit {R}unge--{K}utta methods for {ODEs} with stiff damping term},
  author={Chertock, Alina and Cui, Shumo and Kurganov, Alexander and Wu, Tong},
  journal={SIAM Journal on Numerical Analysis},
  volume={53},
  number={4},
  pages={2008--2029},
  year={2015},
  publisher={SIAM},
doi={10.1137/151005798}
}

@misc{bezanson2015julia,
  title={Julia: A Fresh Approach to Numerical Computing},
  author={Jeff Bezanson and Alan Edelman and Stefan Karpinski and Viral B. Shah},
  year={2015},
  eprint={1411.1607},
  archivePrefix={arXiv},
  primaryClass={cs.MS}
}

@article{rackauckas2017differentialequations,
  title={Differential{E}quations.jl -- {A} performant and feature-rich ecosystem for solving differential equations in {J}ulia},
  author={Rackauckas, Christopher and Nie, Qing},
  journal={Journal of Open Research Software},
  volume={5},
  number={1},
  year={2017},
  publisher={Ubiquity Press},
  doi={10.5334/jors.151},
  url={https://github.com/SciML/OrdinaryDiffEq.jl}
}

@misc{kopecz2023positiveintegrators,
  title={{PositiveIntegrators.jl}: {A} {J}ulia library of positivity-preserving
         time integration methods},
  author={Kopecz, Stefan and Ranocha, Hendrik and contributors},
  year={2023},
  doi={10.5281/zenodo.10868393},
  howpublished={\url{https://github.com/SKopecz/PositiveIntegrators.jl}}
}

@misc{kittisopikul2022interpolations,
  title={{Interpolations.jl}: {F}ast, continuous interpolation of discrete
         datasets in {J}ulia},
  author={Kittisopikul, Mark and Holy, Timothy E and Aschan, Tomas and contributors},
  year={2022},
  howpublished={\url{https://github.com/JuliaMath/Interpolations.jl}}
}

@article{ortleb2017patankar,
  title={Patankar-type {R}unge-{K}utta schemes for linear {PDEs}},
  author={Ortleb, Sigrun and Hundsdorfer, Willem},
  journal={AIP Conference Proceedings},
  volume={1863},
  number={1},
  pages={320008},
  year={2017},
  doi={10.1063/1.4992489}
}

@article{ciallella2022arbitrary,
  title={An arbitrary high order and positivity preserving method for the
         shallow water equations},
  author={Ciallella, Mirco and Micalizzi, Lorenzo and {\"O}ffner, Philipp and
          Torlo, Davide},
  journal={Computers \& Fluids},
  volume={247},
  pages={105630},
  year={2022},
  publisher={Elsevier},
  doi={10.1016/j.compfluid.2022.105630}
}

@misc{ciallella2024high,
  title={A high-order, fully well-balanced, unconditionally 
         positivity-preserving finite volume framework for 
         flood simulations},
  author={Ciallella, Mirco and Micalizzi, Lorenzo and 
          Michel-Dansac, Victor and {\"O}ffner, Philipp 
          and Torlo, Davide},
  year={2024},
  howpublished={\url{https://arxiv.org/abs/2402.12248}}
}

@book{gottlieb2011strong,
  title={Strong stability preserving {R}unge-{K}utta and multistep time
         discretizations},
  author={Gottlieb, Sigal and Ketcheson, David I and Shu, Chi-Wang},
  year={2011},
  publisher={World Scientific},
  address={Singapore}
}

@misc{bartel2024structureRepro,
  title={Reproducibility repository for
         "{S}tructure-Preserving Numerical Methods
           for {F}okker-{P}lanck Equations"},
  author={Bartel, Hanna and Lampert, Joshua and Ranocha, Hendrik},
  year={2024},
  howpublished={\url{https://github.com/JoshuaLampert/2024\_fokker\_planck}},
  doi={10.5281/zenodo.10955446}
}
